\documentclass[12pt]{article}

\usepackage{amssymb,amsmath,theorem}
\usepackage{epsfig}
\usepackage{calc,ifthen}
\usepackage{enumerate}
\usepackage{url}

\newcommand\cF{{\mathcal F}}
\newcommand\cS{{\mathcal S}}

\newcommand\comp[1]{{\mkern2mu\overline{\mkern-2mu#1}}}
\newcommand\rk{\operatorname{rk}}
\newcommand\sbs{\subseteq}

\theoremstyle{change}
{\theorembodyfont{\slshape}
\newtheorem{theorem}{Theorem.}[section]
\newtheorem{lemma}[theorem]{Lemma.}
\newtheorem{corollary}[theorem]{Corollary.}}
\theorembodyfont{\rmfamily}

\def\proof{\noindent{{\sl Proof. }}}

\def\sqr#1#2{{\vbox{\hrule height.#2pt
     \hbox{\vrule width.#2pt height#1pt \kern#1pt
         \vrule width.#2pt}\hrule height.#2pt}}}
\def\eqed{\sqr53}
\def\qed{%
     \ifmmode\eqno\eqed
     \else\nobreak\ \hfill\eqed\medbreak\fi}

\newcommand\Cref[1]{Chapter~\ref{cha:#1}}

\newcommand\lref[1]{Lemma~\ref{lem:#1}}
\newcommand\tref[1]{Theorem~\ref{thm:#1}}
\newcommand\cref[1]{Corollary~\ref{cor:#1}}
\newcommand\sref[1]{Section~\ref{sec:#1}}

\newcommand\pp[1]{\left(#1\right)}

\newcommand\mat[1]{\left(\begin{matrix} #1 \end{matrix}\right)}

\newcommand\1{\mathbf{1}}

\newcommand\wW{\widehat{W}}
\newcommand\wH{\widehat{H}}
\newcommand\wN{\widehat{N}}
\newcommand\wB{\widehat{B}}

\title{Colouring an Orthogonality Graph}

\author{
  C. D. Godsil$\mbox{}^1$\footnotemark[1]~ and M. W. Newman$\mbox{}^2$\footnotemark[1]\footnote{Research supported by NSERC.}\\
  \normalsize $\mbox{}^1$Department of Combinatorics and Optimization\\
  \normalsize University of Waterloo, CANADA\\
  \normalsize $\mbox{}^2$School of Mathematical Sciences\\
  \normalsize Queen Mary, University of London, UK
}

\begin{document}

\maketitle

\abstract{We deal with a graph colouring problem that arises in quantum information theory. Alice and Bob are
each given a $\pm1$-vector of length $k$, and are to respond with $k$ bits. Their responses must be equal if
they are given equal inputs, and distinct if they are given orthogonal inputs; however, they are not allowed
to communicate any information about their inputs. They can always succeed using quantum entanglement, but
their ability to succeed using only classical physics is equivalent to a graph colouring problem. We resolve
the graph colouring problem, thus determining that they can succeed without entanglement exactly when
$k\leq3$.}

\section{Introduction and Background}

We are concerned here with a graph colouring problem that arises in quantum information theory.

The graph $\Omega_{n}$ has vertex set the set of $\pm1$-vectors of length $n$; two vertices are adjacent if they
are orthogonal. Our main result is that the chromatic number of this graph is equal to $n$ if and only if
$n=2^k$ with $k \leq 3$.

This problem arises in the following scenario, introduced in~\cite{BCW98,BCT99}. Alice and Bob are each given
a $01$-vector of length $2^k$, and they are each to respond with a $01$-vector of length $k$. If their input
vectors are equal, then their output vectors must also be equal; if their input vectors differ in exactly
$2^{k-1}$ positions, then their output vectors must be distinct. Also, they are not allowed to communicate any
information about their input vectors to each other.

Their cause is hopeless without some shared resource. Of course this resource must not allow them to share
information about their inputs. We consider two possibilities, which correspond roughly to classical physics
and quantum physics.

If they are allowed to share prior information then they could agree beforehand on a proper colouring of
$\Omega_{2^k}$. They would then each interpret their input as a vertex of $\Omega_{2^k}$, and respond with the
colour of that vertex. Since they are only allowed to output $k$ bits, this only works if $\chi(\Omega_{2^k})
\leq 2^k$.

Now consider that Alice and Bob have some strategy involving some prior shared information, and that their
strategy in guaranteed to succeed. Alice and Bob are then given their respective inputs. Now before they
actually answer, each writes down a list of the response they would have given to all possible inputs. If they
have a winning strategy then they are able to do this. Let $a(x)$ be the entry in Alice's list corresponding
to $x$, and likewise $b(x)$ for Bob. Since Alice and Bob's strategy is guaranteed to succeed for any pair of
inputs, we have $a(x)=b(x)$ for all $x$, and $a(x) \neq b(y)$ whenever $x$ and $y$ correspond to adjacent
vertices. Their answers are restricted to $k$ bits, so their lists contain at most $2^k$ distinct
entries. Thus they have a proper colouring with at most $2^k$ colours. (Note that we have not shown that if
Alice is given the same input on different occasions that she must respond in the same way. Rather, at each
round, her copy of the shared information amounts to a proper colouring.)

In other words, any strategy based on prior shared information is equivalent to colouring $\Omega_{2^k}$ with at
most $2^k$ colours.

If, instead of information, they are allowed a shared resource of \emph{quantum entanglement}, then they can
always succeed when $n=2^k$ for all $k$.  This was first observed by Buhr, Cleve, and Wigderson~\cite{BCW98}
(see also and Brassard, Cleve, and Tapp~\cite{BCT99}).

We do not assume any familiarity with quantum entanglement, qubits and quantum information theory; however,
the interested reader will find a good introduction to these areas in~\cite{NC}. We summarize briefly the
quantum algorithm of~\cite{BCW98,BCT99}. Alice and Bob share between them an entangled quantum state
consisting of $k$ EPR pairs of qubits. The $\pm1$-vectors of length $2^k$ can be interpreted as indexing a
particular family of quantum operations. So given their inputs, they each apply the corresponding operation to
their qubits, and then measure their qubits. They answer the result of their measurements. This turns out to
be a winning strategy. See~\cite{BCW98, BCT99} for a more formal description. Alternatively, see~\cite{mike}
for a version that does not assume any previous background in quantum information theory.

In~\cite{BCW98} it is also shown that for sufficiently large $k$, Alice and Bob cannot succeed by sharing only
prior information. This follows directly from a deep result of Frankl and R\"odl~\cite{FR87} who show that for
large enough $n=4m$ the size of an independent set in $\Omega_{n}$ is at most $(2-\epsilon)^n$ for some
$\epsilon>0$. It follows that the chromatic number must eventually be greater than $n$. We note that the
motivation of~\cite{FR87} has nothing to do with any quantum scenario; furthermore, their result is stronger
than what we state here.

The point of this scenario is then that Alice and Bob can always succeed using quantum physics (i.e., by
sharing quantum entanglement), whereas they cannot always succeed using classical physics (i.e., by sharing
prior information). So our result quantifies the difference between what can be accomplished using quantum or
classical physics for this particular scenario.

The reader may easily verify that $\chi(\Omega_{n}) = n$ for $n=1,2$, and with a little more effort for $n=4$
as well (this case follows trivially from the recursive construction of \sref{recurse}). Unpublished
computations by Gordon Royle determined that $\chi(\Omega_{8}) = 8$ and characterized all of the proper
$8$-colourings. Galliard, Tapp and Wolf~\cite{GTW02} show that the size of an independent set in $\Omega_{16}$
is at most $3912$, which implies its chromatic number is at least $17$.

\section{Some Simple Cases}

It is not hard to see that $\Omega_{n}$ is edgeless if $n$ is odd.

If $n$ is an odd multiple of two, then the vertices can be divided into the \emph{even} vertices (those with
an even number of $-1$'s) and the odd vertices; every edge joins an even vertex to an odd vertex, and so the
graph is bipartite.

If $n=4m$ then every edge joins an even vertex to an even vertex, or an odd vertex to an odd vertex. In fact,
it is not hard to see that these two subgraphs are isomorphic. Furthermore, if a vertex $x$ is adjacent to $y$
then it is also adjacent to $-y$, and $x$ is not adjacent to $-x$. It follows that each component of
$\Omega_{n}$ can be written as a lexicographic product $Y_{n}[\comp{K_2}]$ for some graph $Y$. Note that
$\chi(\Omega_{n}) = \chi(Y_{n})$ and that maximum independent sets in $\Omega_{n}$ are exactly four
copies of maximum independent sets in $Y_{n}$. It follows that $\alpha(\Omega_{n})$ is a multiple of
four. More importantly, it will simplify some of our computational work.

\section{Bounding Independent Sets}\label{sec:boundindsets}

One of the main tools we use in analyzing $\Omega_{n}$ is the Delsarte-Hoffman bound on independent sets
(see~\cite[Section~3.3]{Delsarte} or~\cite[Page 115]{CDS}; alternatively~\cite{mike} for more recent
work).

\begin{theorem}\label{thm:rb}
  Let $X$ be a $d$-regular connected graph on $v$ vertices, and $\tau$ the least eigenvalue of its adjacency matrix. Let $S$
  be an independent set of size $s$ and let $z$ be the characteristic vector of $S$. Then
  \[ s \leq n \frac {-\tau} {d-\tau}. \]
  Furthermore, equality holds if and only if
  \[ A \pp{ z - \tfrac{s}{v}\1} = \tau \pp{ z - \tfrac{s}{v}\1}.\qed\]
\end{theorem}

The graph $\Omega_{n}$ is a graph in the Hamming scheme. We will not go into details here, but the reader is
directed to~\cite[Chapter 12]{blue} for background material on the Hamming scheme and association schemes in
general. The practical consequence of this is that we know a complete set of eigenvectors for
$\Omega_{n}$. For instance, the methods of~\cite[Section12.9]{blue} can be used to establish the following.

\begin{lemma}\label{lem:ham}
  Let $n$ be a multiple of four. Then the least eigenvalue of $\Omega_{n}$ is
  \[ \tau = - \frac{1}{n-1}\binom{n}{\frac{n}{2}} \]
  and the columns of $W$ form a basis for the $\tau$-eigenspace where $W$ is the matrix with rows indexed by
  subsets of $[n]$ and columns indexed by $2$-subsets and $(n\!-\!2)$-subsets, with $(A,p)$-entry
  equal to $(-1)^{|A \cap p|}$.\qed
\end{lemma}

Note that here we are thinking of the vertices of $\Omega_{n}$ as being subsets of $[n]$ instead of
$\pm1$-vectors; two subsets are adjacent when they are at Hamming distance $\frac{n}{2}$.

Let $\wW = \mat{W&\1}$. Putting \tref{rb} and \lref{ham} together, we obtain the following result for our
graphs.

\begin{corollary}\label{cor:rb}
  The size of an independent set in $\Omega_{n}$ is bounded by
  \[ \alpha(\Omega_{n}) \leq \frac{2^n}{n}. \]
  Furthermore, if equality holds then the characteristic vector of a maximum independent set lies in the
  column space of $\wW$.\qed
\end{corollary}

Note that if $\Omega_{n}$ is $n$-colourable, then this bound must hold with equality. In other words, we have
shown that $\chi(\Omega_{n}) > n$ whenever $n=4m$ is not a power of two. It is noteworthy that the question
``When is $\chi(\Omega_{n}) \leq n$?'' becomes trivial when $n$ is not a power of two: the mathematical
analysis is simplest for the cases that are physically uninteresting. Of more immediate use is the fact that
if the bound does not hold with equality, then $\Omega_{n}$ is not $n$-colourable. One way to show that the
bound is not tight is to use the equality condition of \tref{rb}: it suffices to show that there are no
suitable vectors in the $\tau$-eigenspace.

\section{Finding Maximum Independent Sets}\label{sec:findindsets}

Since we know all of the eigenspaces of $\Omega_{n}$, it is not hard to see that $\Omega_{n}$, the even
component of $\Omega_{n}$, and $Y_{n}$ all have the same least eigenvalue $\tau$. If we take only those
columns of $W$ that correspond to the $2$-subsets, and only those rows that correspond to the even vertices,
the resulting column space gives the $\tau$-eigenspace of the even component (this amounts to taking only one
of the two eigenspaces of the Hamming scheme that give the $\tau$-eigenspace on $\Omega_{n}$). If we further
reduce by taking only one vertex (i.e., row) from each pair $\{x,-x\}$, then we obtain a matrix whose columns
form a basis for the $\tau$-eigenspace for $Y_{n}$. We will denote this matrix by $H$. Furthermore, let $\wH =
\mat{H&\1}$.

\begin{corollary}\label{cor:rbY}
  The size of an independent set in $Y_{n}$ is bounded by
  \[ \alpha(Y_{n}) \leq \frac{1}{4}\frac{2^n}{n}. \]
  Furthermore, if equality holds then the characteristic vector of a maximum independent set lies in the
  column space of $\wH$.\qed
\end{corollary}

Let $z$ be the characteristic vector of an independent set $S$ that meets the bound of \cref{rbY}. Then $z =
\wH y$ for some vector $y$. Since $Y_{n}$ is vertex-transitive, we are free to assume that $S$ contains
any particular vertex. As $S$ is maximum, this is equivalent to assuming that $S$ is disjoint from the
neighbourhood of a vertex. Let $\wN$ be the submatrix of $\wH$ with rows corresponding to a neighbourhood;
then we may assume that $z$ takes the value $0$ at the corresponding positions, meaning that $y$ is in the
kernel of $\wN$. Thus we are lead to the following result.

\begin{lemma}
  The kernel of $\wN$ is given by the row space of $\wB=\mat{B&\1}$, where $B$ is the incidence matrix of
  $K_n$.
\end{lemma}

\proof
A direct computation shows that
\[ N B^T = -\1. \]
This means that
\[ \wN \wB^T = 0. \]
We will show that this is the whole kernel by a rank argument.

We note that
\[ BB^T = (n-1)I + J. \]
Thus the eigenvalues of $BB^T$ are $2n-1$ and $n-1$, so $B$ has full rank, and $\rk(B)=\rk(BB^T)=n$. As
$B\1=(n-1)\1$, we see that $\rk(\wB)=n$.

Essentially the same argument determines the rank of $\wN$ as well.  The rows and columns of $N^TN$ are indexed
by $2$-subsets, and the $(a,b)$-entry depends only on $|a \cap b|$. Let $L$ and $\comp{L}$ be the incidence
matrices of the line graph of $K_n$ and its complement, respectively. It follows that
\[ N^TN = c_0 I + c_1 L + c_2 \comp{L}, \]
where
\begin{align*}
  c_0 &= \binom{n}{\frac{n}{2}},\\
  c_1 &= \binom{n}{\frac{n}{2}} - 8\binom{n-3}{\frac{n}{2}-1},\\
  c_2 &= \binom{n}{\frac{n}{2}} - 16\binom{n-4}{\frac{n}{2}-1}.
\end{align*}

The matrices $I,L,\comp{L}$ are simultaneously diagonalizable with known eigenvectors (more precisely: the
line graph of $K_n$ is strongly regular). It follows that the eigenvalues of $N^TN$ are
\[ \frac{n}{2(n-1)}\binom{n}{\frac{n}{2}}, \qquad \frac{n(n-2)}{(n-1)(n-3)}\binom{n}{\frac{n}{2}}, \qquad 0, \]
with respective multiplicities $1$, $\binom{n}{2}-n$, and $n-1$.
As $N\1 = \frac{n}{2}\1$, we see that $\rk(\wN)=\rk(N)$, and so $\dim(\ker(\wN)) = \rk(\wB)$. The result follows.\qed

We note that the rank arguments in the above proof amount to observing that $BB^T$ and $N^TN$ both lie in the
Bose-Mesner algebras of known association schemes.

Let $C$ be the reduced column echelon form of the matrix $\wH \wN$. Then it follows that there are vectors
$x,x'$ such that
\[ z = \wH y = \wH \wB^T x = C x'. \]
Furthermore, since $z$ is a $01$-vector so is $x'$. Since $H$ and $B^T$ have full column rank, it follows that
the rank of $C$ is $n$. Thus it suffices to check all $2^n$ possibilities for $x'$ in order to determine if
there exist any independent sets that meet the bound. We have carried out this computation for $n=8,16$: for
$n=8$, we find that there are eight independent sets of the required size containing a given vertex; for
$n=16$, there are none.

On its own this computation is not particularly satisfying: we have not contradicted Royle's result mentioned
above, and we have established a weaker bound on $\alpha(\Omega_{16})$ than the one given
in~\cite{GTW02}. However, it will turn out that our computations for $n=8,16$ suffice to determine all values
of $n$ for which $\Omega_{n}$ is $n$-colourable. For this purpose, we will rederive the bound of \tref{rb} for
$\Omega_{n}$ twice more.

\section{Colouring $\Omega_{n}$}\label{sec:colour}

It is well-known that for any vertex-transitive graph $X$ on $v$ vertices, $\alpha(X)\omega(X) \leq v$. This
is the \emph{clique-coclique bound}. It also holds for any graph that is a union of classes in an association
scheme. $\Omega_{n}$ falls into both of these categories, but it is also a normal Cayley graph, for which we
can extend this result. In particular, we will show that for a normal Cayley graph $X$, if
$\alpha(X)\omega(X)=|V(X)|$, then $\chi(X)=\omega(X)$. We will need some preliminary results first.

Recall that the \emph{connection set} of a Cayley graph for a group $G$ is the subset $D$ of $G$ such that $a
\sim b$ whenever $ba^{-1} \in D$. If $S$ is a subset of $G$ then we write
\begin{align*}
  S^{-1} &= \{g^{-1} : g \in S\},\\
  Sa &= \{ga : g \in S\}.
\end{align*}

\begin{lemma}
  Let $X$ be a Cayley graph for a group $G$ with connection set $D \sbs G$. Let $S$ be an independent set of
  $X$.
  If $a$ and $b$ are adjacent then $S^{-1}a \cap S^{-1}b = \emptyset$.
\end{lemma}

\proof
Assume $g^{-1}a = h^{-1}b$ for some $g,h \in S$. Then $hg^{-1}=ba^{-1}$. But $h \sim g$ so $hg^{-1} \notin D$,
while $b \sim a$ so $ba^{-1} \in D$.\qed

\begin{corollary}
  If $X$ is a Cayley graph, then $\alpha(X)\omega(X) \leq |V(X)|$.
\end{corollary}

\proof
Let $S$ be an independent set and $C$ a clique. Then by the previous result the sets
\[ S^{-1}c, c \in C \]
are all disjoint.\qed

We note parenthetically that this can be extended to a proof for all vertex-transitive graphs. It is an old
result of Sabidussi~\cite{Sab64} that if $X$ is a vertex-transitive graph, there is an integer $m$ such that
the lexicographic product $X[\comp{K_m}]$ is a Cayley graph. Since
\[ \alpha(X[\comp{K_m}]) = m \alpha(X), \quad \omega(X[\comp{K_m}]) = \omega(X), \quad |V(X[\comp{K_m}])| = m|V(X)|, \]
the result follows for all vertex transitive graphs.

Recall that a Cayley graph is \emph{normal} if its connection set is closed under conjugation. Our purpose in
approaching the clique-coclique bound through Cayley graphs is the following extension from~\cite{cnotes}.
\begin{corollary}\label{cor:normcay}
  If $X$ is a normal Cayley graph and $\alpha(X)\omega(X)=v$, then $\chi(X)=\omega(X)$.
\end{corollary}

\proof
Let $S$ be an independent set and $C$ be a maximum clique. Again, the sets
\[ S^{-1}c, c \in C \]
are disjoint, and so they partition the vertex set. Since $X$ is normal they are also independent sets, and
hence form a colouring.\qed

For our purposes, \cref{normcay} rederives the bound of \tref{rb} for $\Omega_{n}$, but with a different
equality condition. Notice that any clique in $\Omega_{n}$ is a set of pairwise orthogonal vectors, hence
linearly independent, and hence has size at most $n$. Furthermore, an $n$-clique would correspond to a
Hadamard matrix, which certainly exists if $n=2^k$. It follows that we can use \cref{normcay} to rederive the
bound of \tref{rb}, but with a different equality condition.

\begin{corollary}\label{cor:rbc}
  Let $n$ be a power of two.

  The size of an independent set in $\Omega_{n}$ is bounded as
  \[ \alpha(\Omega_{n}) \leq \frac{2^n}{n}. \]
  Furthermore, if equality holds then $\chi(\Omega_{n}) = n$.\qed
\end{corollary}

So $\chi(\Omega_{n})=n$ if and only if $\alpha(\Omega_{n})=\frac{2^n}{n}$. Our determination of
$\alpha(\Omega_{8})$ above is now upgraded to a proof that $\chi(\Omega_{8})=8$.

More generally, we see that not only is it impossible to $n$-colour $\Omega_{n}$ is $n$ is not a power of two,
but if it is possible, then the colouring is exactly a partition of the vertex set into maximum independent
sets that meet the bound of \tref{rb}.

\section{A Recursive Construction}\label{sec:recurse}

The graph $\Omega_{n}$ is an induced subgraph of $\Omega_{2n}$: take exactly those vertices of $\Omega_{2n}$
whose last $n$ entries are the same as the first $n$ entries.
In fact, we can say much more than this.

For vertices $x,r$ of $\Omega_{n}$, let $x^{(r)}$ be the vertex of $\Omega_{2n}$ obtained by concatenating $x$
with the entrywise product of $x$ and $r$ (recall that vertices are $\pm1$-vectors). Let $\Omega_{n}^{(r)}$ be
the subgraph induced by the vertices
\[ \{ x^{(r)} : x \in V(\Omega_{n}) \}. \]
Then we see that $\Omega_{n}^{(r)}$ is isomorphic to $\Omega_{n}$, for any $r$. (The previous example was
$\Omega_{n}^{(\1)}$.) The vertex set of $\Omega_{2n}$ can be partitioned as
\[ \{ V(\Omega_{n}^{r}) : r \in V(\Omega_{n}) \} \]
Furthermore, every vertex of $\Omega_{n}^{(r)}$ is adjacent to every vertex of $\Omega_{n}^{(-r)}$.

Recall that the \emph{join} of two graphs $X_1$ and $X_2$ is $X_1 + X_2 = \comp{\comp{X_1} \cup
\comp{X_2}}$. We have established the following result.

\begin{lemma}
The vertex set of $\Omega_{2n}$ can be partitioned into $2^{n-1}$ copies of $\Omega_{n} + \Omega_{n}$.\qed
\end{lemma}

Note that any independent set in the join of two graphs must lie entirely within one or the other. This gives
us a bound on the size of an independent set in $\Omega_{2n}$: it is at most half of the size of $2^n$ maximum
independent sets in $\Omega_{n}$. We can use this to again rederive the bound of \tref{rb}, with yet another
equality condition.

\begin{corollary}\label{cor:rbr}
  Let $n=2^k$ where $k>1$.

  The size of an independent set in $\Omega_{n}$ is bounded as
  \[ \alpha(\Omega_{n}) \leq \frac{2^n}{n}. \]
  Furthermore, if equality holds then it also holds for $n=2^{k-1}$.\qed
\end{corollary}

We can in fact apply the recursive construction when $n$ is not a power of two. If $n=m2^k$, where $m$ is odd,
then we find that
\begin{equation}\label{m2k}
  \alpha(\Omega_{n}) \leq \frac{2^n}{2^k}.
\end{equation}
When $k>1$ the bound of \cref{rb} is better than \eqref{m2k} by a factor of $m$. When $k=1$, the bound of
\eqref{m2k} is half the number of vertices, and is tight since $\Omega_{n}$ is then bipartite. But for $k=1$
the least eigenvalue of $\Omega_{n}$ is no longer given by \lref{ham}, and applying \tref{rb} in this case
again gives half the number of vertices.

In other words, this recursion does not ever give a tighter bound than \tref{rb}; rather, it is useful because
of the equality condition of \cref{rbr}.

We mention another point of view on this recursion.

For $n=2^k$, let the graphs $\Psi_{n}$ be defined by setting $\Psi_{1}:=\Omega_{1}=\comp{K_2}$, and
recursively defining $\Psi_{2n}$ to be the disjoint union of $2^{n-1}$ copies of $\Psi_{n} + \Psi_{n}$. Then
$\Psi_{n}$ is a spanning subgraph of $\Omega_{n}$, and so \cref{rb} gives a bound on independent sets in
$\Psi_{n}$ as well. It follows that
\begin{align*}
  \alpha(\Psi_{n}) &= \frac{2^n}{n},\\
  \chi(\Psi_{n}) &= n.
\end{align*}
Furthermore we see that $\Psi_{n} = \Omega_{n}$ for $n=1,2,4$, by simply observing that their degrees are the
same. This is an easy way to see that $\chi(\Omega_{4})=4$.

There are of course other edges in $\Omega_{n}$ for general $n$. In fact, $\Psi_{n}$ is not only spanning, it
is asymptotically sparse, in the following sense.

\begin{lemma}
  \[ \lim_{k \to \infty} \frac{|E(\Psi_{2^k})|}{|E(\Omega_{2^k})|} = 0\qed\]
\end{lemma}

In other words, for $n=2^k$, \tref{rb} effectively only ``sees'' the edges of $\Psi_{n}$. This is
expected, since one consequence of the Frankl-R\"odl result is that for large enough $n$, the bound of
\cref{rb} is exponentially too big.

\section{Main Result}

Recall that although we can easily define $\Omega_{n}$ for any positive integer $n$, it is the cases where
$n=2^k$ that are most of interest. In exactly these cases, we have three different ways of proving the same
bound on independent sets: using the Delsarte-Hoffman bound, using the maximum cliques, and a recursive
construction. Furthermore, each approach gives different information in the case where the bound is tight.

We now find that our main result follows directly.

\begin{theorem}
  $\chi(\Omega_{n})=n$ if and only if $n=2^k$ with $k\leq3$.
\end{theorem}

\proof
It follows from \cref{rb} and the comments after it that $\chi(\Omega_{n}) \neq n$ if $n$ is not a power of
two. Furthermore, if $\chi(\Omega_{n}) = n$ then the bound of \cref{rb} holds with equality, and \cref{rbc}
tells us that it is sufficient that this bound holds with equality. Our computations of \sref{findindsets}
deal with the cases $n=2^k$ for $k=3,4$. We then invoke \cref{rbr} to conclude that $\chi(\Omega_{2^k})>2^k$
for all $k>4$.\qed

\section{Further Bounds}

Galliard~\cite{gall} has a construction of an independent set inspired by the methods of Ahlswede and
Khachatrian~\cite{AK97}. We state it in terms of the graph $Y_n$. It is convenient to regard the vertices as
being subsets of $[n]$.

Let $n=2^k$ and $c=\frac{n}{4}-1$. Then the following collection is an independent set in $Y_n$.
\[ \cF_{n} = \{ F \sbs [n] : |F|=2i,\; i \leq 2c ; \; |F \cap [c]| \geq |F \setminus [c]| \} \]
It is not hard to see that $\cF$ is not properly contained in any larger independent set.

It turns out that for $k \leq 3$ this set meets the bound of \cref{rb}, and hence this construction is
maximum. Up to automorphisms of $\Omega_{n}$, this is unique (this follows both from Gordon Royle's
computations and from our work in \sref{findindsets}). Galliard conjectured that $\cF_{n}$ is maximum for all
$n=2^k$.

A recent computation of de Klerk and Pasechnik~\cite{dKP05} using a technique of Schrijver~\cite{Sch04} gives
that $\alpha(Y_{16}) \leq 576$, which is exactly the size of $\cF_{16}$. The reader is referred
to~\cite{dKP05} for details and further results. We do not know of any determination of $\alpha(\Omega_{4m})$
for $m>4$.

There is another way to look at the the collection $\cF$. If we replace each element $F$ of $\cF_{n}$ with its
symmetric difference with $[c]$ we obtain the collection of all odd subsets of $[n]$ of size at most $c$. More
generally, if we assume only that $n=4m$, we have the following construction.
\[ \cS_{n} = \{ F \sbs [n] : |F| \not\equiv m \!\!\!\!\pmod 2; \; |F| < m \} \]
We conjecture that these are in fact maximum in general. If this is true, it would imply our statement of the
Frankl-R\"odl result, namely that
\[ \alpha(\Omega)_{4m} \leq (2-c)^{4m} \]
for some $c>0$, for large enough $m$.

\bibliographystyle{plain}
\bibliography{ortho}

\begin{thebibliography}{10}

\bibitem{AK97}
R.~Ahlswede and L.~H. Khachatrian.
\newblock The complete intersection theorem for systems of finite sets.
\newblock {\em European J. Combin.}, 18(2):125--136, 1997.

\bibitem{BCT99}
G.~Brassard, R.~Cleve, and A.~Tapp.
\newblock Cost of exactly simulating quantum entanglement with classical
  communication.
\newblock {\em Phys. Rev. Lett.}, 83(9):1874--1877, 1999.

\bibitem{BCW98}
H.~Buhrman, R.~Cleve, and A.~Widgerson.
\newblock Quantum vs. classical communication and computation.
\newblock In {\em Proceedings of the 30th Annual {ACM} Symposium on the Theory
  of Computing}, pages 63--68, 1998.

\bibitem{CDS}
D.~M. Cvetkovi{\'c}, M.~Doob, and H.~Sachs.
\newblock {\em Spectra of graphs}.
\newblock Academic Press Inc. [Harcourt Brace Jovanovich Publishers], New York,
  1980.

\bibitem{dKP05}
E.~de~Klerk and D.~V. Pasechnik.
\newblock A note on the stability number of an orthogonality graph.
\newblock 2004.
\newblock Available as \mbox{\texttt{arXiv:math.CO/0505038}}.

\bibitem{Delsarte}
P.~Delsarte.
\newblock An algebraic approach to the association schemes of coding theory.
\newblock {\em Philips Res. Rep. Suppl.}, (10):vi+97, 1973.

\bibitem{FR87}
Peter Frankl and Vojt{\v{e}}ch R{\"o}dl.
\newblock Forbidden intersections.
\newblock {\em Trans. Amer. Math. Soc.}, 300(1):259--286, 1987.

\bibitem{gall}
V.~Galliard.
\newblock Classical pseudo-telepathy and colouring graphs.
\newblock Master's thesis, ETH Zurich, 2001.

\bibitem{GTW02}
V.~Galliard, A.~Tapp, and S.~Wolf.
\newblock The impossibility of pseudo-telepathy without quantum entanglement.
\newblock 2002.
\newblock Available as \mbox{\texttt{arXiv:quant-ph/0211011}}.

\bibitem{blue}
C.~D. Godsil.
\newblock {\em Algebraic Combinatorics}.
\newblock Chapman \& Hall, New York, 1993.

\bibitem{cnotes}
C.~D. Godsil.
\newblock {Interesting Graphs and Their Colourings}.
\newblock Unpublished notes, 2003.

\bibitem{mike}
M.~W. Newman.
\newblock {\em Independent Sets and Eigenspaces}.
\newblock PhD thesis, University of Waterloo, 2004.

\bibitem{NC}
M.~A. Nielsen and I.~L. Chuang.
\newblock {\em Quantum computation and quantum information}.
\newblock Cambridge University Press, Cambridge, 2000.

\bibitem{Sab64}
G.~Sabidussi.
\newblock Vertex-transitive graphs.
\newblock {\em Monatsh. Math.}, 68:426--438, 1964.

\bibitem{Sch04}
A.~Schrijver.
\newblock New code upper bounds from the terwilliger algebra and semidefinite
  programming.
\newblock 2004.
\newblock Preprint, available at \mbox{\url{http://homepages.cwi.nl/~lex/}}.

\end{thebibliography}

\end{document}